\input amstex
\documentstyle{amsppt}
\document
\topmatter
\title
K\"ahler manifolds with some Killing tensors.
\endtitle
\author
W\l odzimierz Jelonek
\endauthor

\abstract{The aim of this paper is to classify compact, simply
connected K\"ahler manifolds which admit J-invariant Killing tensor with two eigenvalues of multiplicity 2 and n-2 and with constant eigenvalue corresponding to 2-dimensional eigendistribution.}
\thanks{MS Classification: 53C55,53C25. Key words and phrases:
K\"ahler manifold, holomorphic foliation, homothetic foliation,
special K\"ahler-Ricci potential, special K\"ahler
potential}\endthanks
 \endabstract
\endtopmatter
\define\G{\Gamma}
\define\DE{\Cal D^{\perp}}

\define\n{\nabla}
\define\om{\omega}
\define\r{\rightarrow}
\define\w{\wedge}
\define\k{\diamondsuit}
\define\th{\theta}

\define\a{\alpha}

\define\lb{\lambda}

\define\1{D_{\lb}}
\define\2{D_{\mu}}
\define\0{\Omega}

\define\De{\Cal D}

\define\m{(M,g,J)}

\bigskip
{\bf 0. Introduction. } The aim of the present paper is to
classify compact K\"ahler manifolds $(M,g,J)$,
$dim M=2n>2$,  admitting a $J$-invariant Killing tensor $S$ with two eigenvalues $\lb,\mu$
of multiplicity $2$ and $2n-2$ and with a constant eigenvalue $\lb$ corresponding to the $2$-dimensional eigendistribution.

We show that if $\mu$ is a nonconstant function then such manifold admit
 a complex,
holomorphic, totally geodesic foliation $\Cal F$ by curves on an open dense subset of $M$ (corresponding to the eigendistribution $\1$ of $S$) and is a $\Bbb {CP}^1$-holomorphic bundle $\Bbb P(L\oplus\Cal O)$ over a K\"ahler Hodge manifold $(N,g,J)$ or is a projective space
$\Bbb{CP}^n$.  In the first case both eigenvalues of $S$ are everywhere distinct and in the second case they coincide in an exactly one point.

If both eigenvalues of $S$ are constant then the totally geodesic distribution corresponding to $\1$ may not be holomorphic.
If $dim M=4$ then also in the case of constant eigenvalues $\Cal F$ is holomorphic and in that case if $\mu$ is constant then the simply connected covering space $\tilde M$ of $M$ is the product   $\Sigma\times N$  and the lift of $S$ to $\tilde M$ has eigendistributions $T\Sigma,TN$. If $S$ is a Killing tensor with eigenvalues $\lb,\mu$ where $\lb$ is constant then  $S'=S-\lb I$ is a Killing tensor with eigenvalues $0,\mu'$.  Hence we can assume that $\lb=0$.  If $S$ is a J-invariant Killing tensor with eigenvalues $0,\mu$ of multiplicity $2, 2n-2$ respectively  then it is known that $\phi(X,Y)=S(JX,Y)-\mu\om(X,Y)$ is a Hamiltonian form if  the form  $S(JX,Y)-2\mu\om(X,Y)$ is closed  (see [A-C-G-1], p.407). This last condition is satisfied if $dim M=4$.  In fact in [A-C-G-2] it is proved that  if $dim M=4$ and $S$ is a symmetric $J$ invariant tensor with arbitrary eigenvalues $\lb,\mu$ then  $\phi(X,Y)=S(JX,Y)-(\mu+\lb)\om(X,Y)$ is a Hamiltonian form if  and only if $S$ is a  Killing tensor.
 We shall prove in the present paper that if $\mu$ is non-constant eigenvalue of a Killing tensor with eigenvalues $0,\mu$ of multiplicity $2,2n-2$ respectively then the tensor $\phi(X,Y)=S(JX,Y)-\mu\om(X,Y)$ is a Hamiltonian K\"ahler 2-form also if dim $M>4$.

It is true also for constant $\mu$ if $dim M=4$. However if $dim M>4$ we do not know whether this relation holds if $\mu$ is constant. Note that Hamiltonian forms are classified.
\bigskip
{\bf 1. Killing tensors.}  Let us recall that a (1,1) symmetric  tensor $S$ on a Riemannian
manifold $(M,g)$ is called a Killing tensor if $$g(\n S(X,Y),Z)+g(\n S(Z,X),Y)+g(\n S(Y,Z),X)=0$$ for all $X,Y,Z\in TM$.  Equivalently it means that  $g(\n S(X,X),X)=0$  for all $X\in TM$.
 Define  the integer-valued function  $E_S(x)=($ the number of distinct eigenvalues of $S_x$) and
 set $M_S=\{ x\in M:E_S$ is constant in a neighborhood of $x\}$.
 The set $M_S$ is open and dense in $M$ and the eigenvalues $\lb_i$ of $S$ are distinct and
 smooth in each component $U$ of $M_S$.  Let us denote by $D_{\lb_i}$
 the eigendistributions corresponding to $\lb_i$. We have  (see [J-3])
\medskip
{\bf Proposition 1.1.} {\it Let } $S$ {\it be a Killing  tensor on} $M$ {\it and U be a component of }
 $M_S$ {\it and} $\lb_1,\lb_2,..,\lb_k \in C^{\infty}(U)$
{\it be eigenfunctions of } $S$. {\it Then for all} $X\in D_{\lb_i}$ {\it we have}
   $$ \n S(X,X)=-\frac12 \n \lb_i\parallel X \parallel ^2\tag 1.1$$
{\it and}  $D_{\lb_i}\subset \ker d\lb_i$. { \it If } $i\ne j$ {\it and} $X\in \G(D_{\lb_i}),Y \in \G(D_{\lb_j})$  {\it then}
$$ g(\n_X X ,Y)=\frac12\frac{Y\lb_i}{\lb_j-\lb_i}\parallel X \parallel ^2.\tag 1.2$$

If $T(X,Y)=g(SX,Y)$ is a Killing tensor on $(M,g)$ and $c$ is a geodesic on $M$ then the function  $\phi(t)=T(\dot c(t),\dot c(t))$ is constant on the domain of $c$. In fact $\phi'(t)=\n_{\dot c(t)}T((\dot c(t),\dot c(t))=0$.

A distribution $\Cal D\subset TM$ we call umbilical if $2\n_XX_{|\DE}=g(X,X)\xi$ for every $X\in \G(\De)$ and certain $\xi\in\G(\DE)$. A distribution $\De$ is called totally geodesic if is umbilical with $\xi=0$.

{\bf  Lemma 1.}  {\it  Let $\De$ be a complex, umbilical distribution $J\De\subset \De$ and $dim\De=2$.  Then
$$2\n_XY_{|\DE}=g(X,Y)\xi+\om(X,Y)J\xi$$ for all $X,Y\in\G(\De)$.}

\medskip
{\it Proof.} It is clear that  $(\n_XY+\n_YX)_{|\DE}=g(X,Y)\xi $.  Since  $2\n_{JX}JX_{|\DE}=g(X,X)\xi$ it follows that $2\n_{JX}X_{|\DE}=-g(X,X)J\xi$.
Hence $$[X,JX]_{|\DE}=g(X,X)J\xi$$ and $[X,Y]_{|\DE}=\om(X,Y)J\xi$  for all $X,Y\in \G(\De)$. $\k$

Thus we obtain

{\bf Lemma 2.} {\it Let $\De$ be a complex, totally geodesic distribution  (i.e.  $\n_XX_{|\DE}=0$ for all $X\in\G(\De)$) and $dim\De=2$. Then $\De$ is a totally geodesic foliation.}

\medskip
{\it Proof.}  It follows from Lemma 1, in our case $\xi=0$.$\k$
\medskip
{\it Corollary.}  Let $S$ be a $J$-invariant Killing tensor $S$ with two eigenvalues $\lb,\mu$
of multiplicity 2 and n-2 and with a constant eigenvalue $\lb$ corresponding to the 2 - dimensional eigendistribution $\De$ on a K\"ahler manifold $(M,g,J)$.  Then   $\De$ is integrable and totally geodesic complex distribution.
\medskip
{\it Remark.} We can always assume that $\lb=0$ considering the Killing tensor  $S'(X,Y)=S(X,Y)-\lb g(X,Y)$.

A hamiltonian 2-form is a $J$ invariant 2-form on the K\"ahler manifold $(M,g,J)$ such that there exists a function $\sigma$ on $M$ such that
$$\n_X\phi=\frac12(d\sigma\w JX-d^c\sigma\w X).\tag 1.2$$  It is easy to check that  $d\sigma=d tr_{\om}\phi$ and we shall assume that $\sigma=tr_{\om}\phi$.  If $\phi$ is a hamiltonian form then the tensor $S(X,Y)=\phi(X,JY)-\sigma g(X,Y)$ is a Killing tensor.  In fact from  (1.2) it easily follows that $\frak C_{X,Y,Z}\n_X\phi(Y,JZ)=\frak C_{X,Y,Z}d\sigma(X)g(Y,Z)$  where $\frak C_{X,Y,Z}$  means cyclic sum.   Hence if   $\phi$ is a Hamiltonian form  such that the symmetric tensor $\phi(X,JY)$ has two eigenvalues  $\lb,0$  of multiplicity  $2,n-2$ respectively  with $\mu=0$ constant  then the corresponding Killing tensor $S(X,Y)=\phi(X,JY)-\sigma g(X,Y)$  has eigenvalues $0,-\lb$ corresponding to the same eigendistributions.

\medskip
{\bf 2. Complex homothetic foliations.} We start with  (see [V], [J-1]
):

{\bf Definition.} A foliation $\Cal F$ on a Riemannian manifold
$(M,g)$ is called conformal if
$$L_Vg=\th(V)g$$ holds on $T\Cal F^{\perp}$ where $V\in\G(T\Cal F)$ and $\th$ is a one
form vanishing on $T\Cal F^{\perp}$.  A foliation $\Cal F$ is
called homothetic if it is conformal and $d\th=0$.
\medskip
{\bf Definition.} A foliation $\Cal F$ on a Riemannian manifold
$(M,g)$ is called holomorphic if $L_VJ(TM)\subset T\Cal F$ for every $V\in\G(T\Cal F)$.
\medskip
{\bf Theorem 2.1}   {\it Let us assume that $\Cal F$ is a conformal, homothetic , totally geodesic foliation  by curves.    Then $\Cal F$ is holomorphic in $U=\{x\in M:\th_x\ne0\}$.}

\medskip
{\it Proof.}  We have $2\n_XX=-g(X,X)\th^{\sharp}$ for $X\in T\Cal F^{\perp}$.  On the other hand for $X,Y\in\G(T\Cal F^{\perp})$  $\th([X,Y])=-d\th(X,Y)=0$ hence   $[X,Y]_{|T\Cal F}=\a(X,Y)J\th^{\sharp}$ for some two form $\a$ on $T\Cal F^{\perp}$.  Hence
$2\n_XY=-g(X,Y)\th^{\sharp}+\a(X,Y)J\th^{\sharp}$ and $2\n_XJY=-g(X,JY)\th^{\sharp}+\a(X,JY)J\th^{\sharp}$.  On the other hand $$2\n_XJY=-g(X,Y)J\th^{\sharp}-\a(X,Y)\th^{\sharp}.$$  Hence   $\a(X,Y)=g(X,JY)=-\om(X,Y)$ and
$2\n_XY=-g(X,Y)\th^{\sharp}-\om(X,Y)J\th^{\sharp}$.  Let us recall that $\Cal F$ is holomorphic if  $$\n_{JX}V-J\n_XV\in\G(\De)\tag 2.1$$ for every $X\in TM$ and $V\in\G(\De)$ where $\De=T\Cal F$. If $X\in\G(\De)$ then (2.1) is satisfied since $\Cal F$ is totally geodesic.
Let $X,Y\in\G(\De^{\perp})$.  We have to show that  $g(\n_{JX}V-J\n_XV,Y)=0$ or equivalently that $g(V,\n_{JX}Y+\n_XJY)=0$.  But we have
$\n_{JX}Y+\n_XJY_{|\De}=-g(X,JY)\th-\om(JX,Y)J\th-g(X,JY)\th-\om(X,JY)J\th=0$ which means that $\Cal F$ is holomorphic.$\k$
\medskip
{\it Remark.}   Note that if $dim M=4$ then $\Cal F$ is also holomorphic when $\th=0$.  It follows from Lemma 2 that $\De^{\perp}$ is then totally geodesic and integrable hence $(\n_{JX}Y+\n_XJY)_{|\De}=0$ for $X,Y\in\G(\De^{\perp})$. It follows that the following theorem holds
\medskip
{\bf Theorem 2.2} {\it Let $\m$ be a K\"ahler manifold, dim M=4 and let $\Cal F$ be a complex, conformal, homothetic and totally geodesic foliation on $M$. Then $\Cal F$ is holomorphic.}

\medskip
{\bf Theorem 2.3} {\it Let $\m$ be a compact K\"ahler manifold and let $\Cal F$ be a complex, conformal, homothetic and totally geodesic foliation on $M$ such that $\th$ does not vanish identically on $M$. Then $\Cal F$ is holomorphic.}

\medskip
{\it Proof.} Let $U=\{x\in M:|\th_x|\ne0\}$.  Then exactly as in [J-1], [J-2] we show that  $U$ is dense in $M$.  Since $\Cal F$ is holomorphic on an open and dense subset it follows that it is holomorphic in the whole of $M$.$\k$
\bigskip
{\bf 3.} {\bf Examples of Killing tensor fields on compact K\"ahler manifolds.} First we give a definition
 \medskip {\bf Definition.} {\it A nonconstant
function $\tau\in C^{\infty}(M)$, where $\m$ is a K\"ahler
manifold, is called a special K\"ahler potential if the field $X=
J(\n \tau)$ is a Killing vector field and, at every point with
$d\tau\ne 0$ all nonzero tangent vectors orthogonal to the fields
$X,JX$ are  eigenvectors of  $\n d\tau$. If  $dim M=4$ we shall additionally assume that the field $\n\tau$ is an eigenfield of Ricci tensor $Ric$ of $(M,g,J)$.}
\medskip
  K\"ahler compact manifolds which admit special K\"ahler potential are in fact classified in [D-M-2]  ( see [J-1]). These are $\Bbb{CP}^1$ holomorphic bundles $\Bbb P(L\oplus\Cal O)$ over Hodge K\"ahler manifolds, where $L$ is a holomorphic line bundle over K\"ahler Hodge manifold $(N,\om,J)$  or $\Bbb{CP}^n$.

Let $\tau$ be a special K\"ahler potential on a compact  K\"ahler manifold
$\m$. Then (see [D-M-1], Prop.11.5) $\tau$ has exactly two critical manifolds, which are the $\tau$-preimages of its extremum values $\tau_{max},\tau_{min}$ and one of the following two cases holds:

(a)  both critical manifolds of $\tau$ are of complex codimension one;

(b) one of the critical manifolds has complex  codimension $1$ and the other is a single point.

What is more the function $Q=g(\n\tau,\n\tau)$ is a composite consisting of $\tau$ followed by a $C^{\infty}$ function $[\tau_{min},\tau_{max}]\ni \tau \r Q\in \Bbb R$ which is positive on the open interval $(\tau_{min},\tau_{max})$ and vanishes at the endpoints while the values of $\frac{dQ}{d\tau}$ at the endpoints are mutually opposite and nonzero. We shall denote $|\frac{dQ}{d\tau}|=2a\ne0$.
Let us denote $$\Cal V=span\{\n\tau,J\n\tau\}$$ on $U=\{x\in M:d\tau(x)\ne0\}$ and
let $\Cal F$ be a foliation on $U$ given by the integrable
distribution $\Cal V$.  By $\Theta$ we denote the
eigenvalue of the Hessian $H^{\tau}$ corresponding to the
distribution $\Cal H=\Cal V^{\perp}$. Then $\Cal F$ is a totally geodesic,
holomorphic complex conformal foliation. We have $\th =2\frac \Theta Q
d\tau,\zeta=2\frac \Theta Q \n\tau$, where $\th(X)=g(\zeta,X)$ and $Q=g(\n\tau,\n\tau)$.  There exists a constant $c$ such that $$\frac Q{\Theta}=2(\tau-c)\tag 3.1$$ or $\Theta=0$.  We shall assume that  $\Theta\ne0$ on $U$.  Note also that  $$ dQ=2\Lambda d\tau.\tag 3.2 $$
\medskip
Let us assume that $\Theta\ne0$ and  define a tensor $S$ on $M$ by  $SX=0$ if $X\in \Cal V$  and $SX=(\tau-c) X$ if $X\in \Cal H$. If $M=\Bbb P(L\oplus\Cal O)$ then both distributions  $\Cal V,\Cal H$ extends to smooth distributions over $M$   and clearly $S$ is a smooth tensor on $M$.  In that case the foliation tangent to $\Cal V$ is a conformal, homothetic foliations by curves and  $\th=d\ln|\tau-c|$.
If  $M=\Bbb{CP}^n$ at the point  $x_0$  where $\tau(x_0)=c$ the distributions $\Cal H,\Cal V$ are not defined and we define $S_{x_0}=0$.  If $(M,g,J)$ with special Killing potential $\tau$ admits a critical submanifold consisting of one point $x_0$  then $\tau(x_0)=c$ and
$$g(SJX,Y)=(\tau-c)\om-\frac{\tau-c}Q d\tau\w d^c\tau$$  which means that $S$ is smooth also in a neighborhood of $x_0$ since $\lim_{x\r x_0}\frac{\tau-c}Q=\frac1{\frac{dQ}{d\tau}}=\pm\frac1{2a}$ is finite and $\frac{\tau-c}Q$ is smooth in a certain neighborhood of $x_0$.
Note that at points where $d\tau(x)=0$ but $\tau(x)\ne c$  the distributions $\Cal H,\Cal V$ are defined.

We also define Killing tensors with two constant eigenvalues on the products $M=\Sigma\times N$ where $\Sigma$ is a complex surface
and $N$ is a K\"ahler manifold with $dim N=n-2$  just defining   $\Cal V=T\Sigma,\Cal H=TN$  and  $SX=0$ for $X\in \Cal V$ and $SX=\mu X$ for $X\in \Cal H$ where $\mu\in\Bbb R$ is any real number.

We prove now that the killing tensors constructed above are related to Hamiltonian forms  namely the 2-form  $\phi=(\tau-c)\frac {d\tau\w d^c\tau}Q$  is a Hamiltonian form  (see also  [A-C-G] p.368 )  which  means that (note that $\sigma=tr_{\om}\phi=\tau-c$)
$$\n_X\phi =\frac12(d\tau\w JX-d^c\tau\w X).\tag 3.4$$
Note that  $\phi$  is defined on $U=\{x\in M:d\tau(x)\ne0\}$ and has a smooth extension on the whole of $M$.  On the critical manifold $N$ with $dim_{\Bbb C} N=n-1$ it is clear since this form is $(\tau-c)\om_1$ where $\om_1$ is the volume form of the foliation $\Cal F$ and $\Cal F$  extends smoothly to $N$.   If the critical manifold is $\{x_0\}$  then the function  $\frac{\tau-c}Q$ is smooth at $x_0$ and hence $\phi$ extends smoothly to $x_0$.
  We have

$$ \gather \n_X\phi= d\tau(X)\frac {d\tau\w d^c\tau}Q -\frac{dQ(X)}{Q^2}(\tau-c)d\tau\w d^c\tau\tag 3.5\\+\frac{\tau-c}Q\n_Xd\tau\w d^c\tau+\frac{\tau-c}Qd\tau\w \n_X d^c\tau. \endgather $$
We consider  two cases:

(a)If $X\in \De^{\perp}$ then (note that $\frac Q{\tau-c}=2\Theta$)
$$\gather \n_X\phi(Y,Z)=\frac1{2\Theta}(\Theta g(X,Y)d^c\tau(Z)-\Theta g(Z,X)d^c\tau(Y))+\\ \frac1{2\Theta}(\Theta g(JX,Z)d\tau(Y)-\Theta g(JX,Y)d\tau(Z))\endgather$$  hence  $\n_X\phi(Y,Z) =\frac12(d\tau\w JX-d^c\tau\w X)(Y,Z).$

(b) If $X\in \De$  then we can consider two subcases:

(1) $X=J\n \tau$  In this case  one can easily check that  $\n_X\phi=0$ hence again (3.4) holds,

(2) $X=\n\tau$  Then
$$\gather \n_X\phi(Y,Z)=  d\tau\w d^c\tau(Y,Z)- \frac{2\Lambda}{2\Theta}d\tau\w d^c\tau(Y,Z)\\ +\frac{\Lambda}{2\Theta}(d\tau(Y)d^c\tau(Z)-d\tau(Z)d^c\tau(Y))+\frac{\Lambda}{2\Theta}(d\tau(Y)d^c\tau(Z)-d\tau(Z)d^c\tau(Y))\\=d\tau\w d^c\tau(Y,Z)\endgather$$
 and again (3.4) holds.
\bigskip
{\bf 4. Classification of   Killing tensors.}  Let $(M,g,J)$ be a compact K\"ahler manifold of real dimension $2n$. Let $S$ be a J-invariant Killing tensor on $(M,g,J)$  with two eigenvalues $0,\mu$ and let $U=\{x\in M:\mu(x)\ne0\}$. We assume that in $U$ $dim ker S=2, dimker(S-\mu I)=2n-2$. Then  $\De=ker S$ is a totally geodesic complex foliation defined in $U$ which is holomorphic in the set $V=\{x\in U:\n \mu(x)\ne 0\}$ in view of Th.2.1.  In fact since in $U$ we have  (see (1.2)) $\n_XX_{|\2}=0$ if $X\in\G(\De)$ and $\n_XX_{|\De}=-g(X,X)\frac{\n\mu}{2\mu}$ if $X\in\G(\2)$ it follows that in $U$ $\De$ is totally geodesic, conformal and homothetic foliation with $\th=d\ln|\mu|$.  Let  $K=\{x\in M:\mu(x)=0\}$.

Just as in [J-1] we prove that in $V=\{x\in U:d\mu(x)\ne0\}$ there exist locally defined holomorphic Killing vector field with special K\"ahler potential $\tau$  such that $\th=d\ln|\tau-c|=d\ln|\mu|$. We can assume that $\tau-c=\mu$.  Hence we see that the field $\xi=J(\n\mu)$ is a holomorphic Killing vector field defined in the whole of $V$.
 \medskip
 {\bf  Proposition 4.1.}  {The set  $V$ is connected and the set  $M-V$ has an empty interior.}
 \medskip
 {\it Proof.}  Note that the function $\mu$  is the special Killing potential in $V$. Note that   $\frac Q{\Theta}=2\mu$ where $Q=g(\xi,\xi)$  and $\Theta$ is an eigenvalue of $H^{\mu}$ corresponding to the distribution $\De^{\perp}$. The function $\Theta$ is bounded on $M$ and   $|\Theta|\le sup_{X\in SM}|H^{\mu}(X,X)|$ where $SM=\{X\in TM:g(X,X)=1\}$ is compact unit subbundle of $TM$.  Next we show that the interior of the set $M-V$ is empty and the set $M-V$ is connected.  Let $V_1$ be a connected component of $V$ and let $x_0\in V_1$.   Let  $x\in (M-V)\cup (V-V_1)$.  Let $c$ be a geodesic
joining the points $c(0)=x_0$ and $c(k)=x_1$.   Let $J$ be the Jacobi field along $c$ which in $V$ coincides with $\xi$.    We show that $g(J,\dot{c})=0$. Note that $c$ intersects the set   $K=\{x\in M:\mu(x)=0\}$ or $\{x\in U:|\th|=0\}$.  Let $l\in dom c$ be such a number that $c([0,l))\subset V$  and  $c(l)\in K$.  Since  $Q(c(t))=\mu\Theta(c(t))$ and $|\Theta|$ is bounded on $M$ it follows that $\lim_{t\r l-}Q(c(t))=0$.   On the other hand   if  $c(l)\in \{x\in U:|\th|=0\}$ and  $c([0,l))\subset V$ then exactly as in [J-1], [J-2]]   we prove that $\lim_{t\r l-}Q(c(t))=0$.  Since  $g(J,\dot c)$ is constant along $c$ it follows that  $g(J,\dot c)=0$.  Now let us assume that there exists an open neighborhood $W$ of $x$ such that $W\subset  (M-V)\cup (V-V_1)$. Let us assume that  $x=\exp_{x_0}X_0$  for some $X_0\in T_{x_0}M$.  Then  $g(X_0,\xi(x_0))=0$.   We find an open neighborhood $W'$ of $X_0$ such that $\exp_{x_0}(W')\subset W$.  We find $X\in W'$ such that $g(\xi(x_0),X)\ne 0$.  Hence the geodesic $d(t)=\exp_{x_0}tX$ intersects the boundary of $V$ and   $g(J,\dot d)\ne 0$ a contradiction.   It follows that $int  (M-V)\cup (V-V_1)=\emptyset$ which means that $V=V_1$ and $int(M-V)=\emptyset$. $\k$

\medskip

{\bf Proposition 4.2.}  { The field  $\xi=J\n\mu$ is a holomorphic Killing vector field on the whole of $M$ with a special Killing potential $\mu$.}

\medskip
{\it Proof.}  Note that  $H^{\mu}(JX,JY)=H^{\mu}(X,Y)$ on the open and dense subset $V$ of $M$.    From the continuity of $H^{\mu}$ it follows that   $H^{\mu}(JX,JY)=H^{\mu}(X,Y)$ in the whole of $M$ which proves that $\xi=J(\n\mu)$ is a holomorphic Killing vector field on $M$.  Since on $M-V$  the field $\n \mu$ vanishes and $\mu$ is a special Killing potential on $V$ it follows that $\mu$ is a special Killing potential on $M$. $\k$

{\bf Theorem  4.1.} {\it Let assume that  a Killing tensor $S$ on a compact K\"ahler manifold $(M,g,J)$ has two eigenvalues $0,\mu$ of multiplicity $2,2n-2$ respectively, where $2n=dim M$. Let us assume that $\mu$ is not constant. Then $\mu\ne 0$ on the whole of $M$ and $M=\Bbb P(L\oplus \Cal O)$ where  $L$ is a holomorphic line bundle over Hodge K\"ahler manifold $(N,h,J)$ with Chern connection with curvature  $\Omega=2a\om_h$   or there exists one point $x_0\in M$ such that $\mu(x_0)=0$ and then $M$ is biholomorphic to the projective space  $\Bbb{CP}^n$.  In both cases manifold $M$ admits special K\"ahler potential $\mu$ and $\phi=S(JX,Y)-\mu\om(X,Y)$ is a global Hamiltonian 2-form such that $\phi=\mu\frac{d\mu\w d^c\mu}Q$  on $U=\{x:d\mu(x)\ne0\}$. }

\medskip
{\it Proof. }
In the papers [D-M-1], [D-M-2] in fact there is given a classification of compact K\"ahler manifolds admitting a special K\"ahler potential if $dim M\ge 6$ (see [J-1]). Note also that if $dim M=4$   function $\mu$ is a special K\"ahler- Ricci potential (see [D-M-2],[J-1]).  Hence it follows that the Killing vector field $J\n \mu$ has exactly two critical submanifolds which are both of (complex) codimension 1  and then is biholomorphic to the manifold  $\Bbb{P}(L\oplus\Cal O)$ where $L$ is a holomorphic line bundle over Hodge manifold $(N,\om_N)$   with Chern connection with curvature  $\Omega=2a\om_h$  or one of two critical manifold is a point $\{x_0\}$ and then   $M=\Bbb{CP}^n$. Note that the point $x_0$ we have $\mu(x_0)=0$ and this is the only point with this property  ( the special K\"ahler potential  $\tau=\mu$  satisfies $\tau(x_0)=c=0$ at $x_0$).

If dim $M=4$ and  $\mu$ is constant then simply connected covering space $\tilde M$ of $M$ is a product $\Bbb{CP}^1\times N$ since in that case the foliation $\Cal F$ is holomorphic  (see Th.2.2).$\k$

{\bf Theorem 4.2.} {\it Let assume that  a Hamiltonian 2-form  $\phi$ on a compact K\"ahler manifold $(M,g,J)$ has two eigenvalues $\mu,0$ of multiplicity $2,2n-2$ respectively, where $2n=dim M$. Let us assume that $\mu$ is not constant. Then $\mu\ne 0$ on the whole of $M$ and $M=\Bbb P(L\oplus \Cal O)$ where  $L$ is a holomorphic line bundle over Hodge K\"ahler manifold $(N,h,J)$ with Chern connection with curvature  $\Omega=2a\om_h$   or there exists one point $x_0\in M$ such that $\mu(x_0)=0$ and then $M$ is biholomorphic to the projective space  $\Bbb{CP}^n$.  In both cases manifold $M$ admits the special K\"ahler potential $\tau$ such that a Hamiltonian 2-form  $\phi=(\tau-c)\frac{d\tau\w d^c\tau}Q$  on $U=\{x:d\tau(x)\ne0\}$.  If $\mu$ is constant then the simply connected covering space $\tilde M$ of $M$  is the product $\tilde M= \Sigma\times N$ and the lift of $\phi$ to $M$ has eigendistributions  $T\Sigma,TN$.}

\medskip
{\it Remark.} If a compact K\"ahler manifold  $(M,g,J)$ of dimension $2n\ge4$  admits a special K\"ahler potential $\tau$ such that if $dim M=4$ the field  $\n\tau$ is an eigenfield of the Ricci tensor $Ric$ of $(M,g,J)$ and the eigenvalue $\Theta$ of $H^{\tau}$ corresponding to the eigendistribution $\Cal H$ is nonzero then the 2-form $\phi=(\tau-c)\frac{d\tau\w d^c\tau}Q$ defined on $U=\{x:d\tau(x)\ne0\}$ extends to the smooth Hamiltonian 2-form on the whole of $M$.

\bigskip
\centerline{\bf References.}

\par
\medskip

\cite{A-C-G-1}  V. Apostolov, D.M.J. Calderbank, P. Gauduchon  {\it Hamiltonian 2-forms in K\"ahler geometry, I General Theory} J. Diff. Geom 79, (2006), 359-412
\par
\medskip

\cite{A-C-G-2}  V. Apostolov, D.M.J. Calderbank, P. Gauduchon {\it Ambitoric geometry I: Einstein metrics and extremal ambiK\"ahler structures}   	J. reine angew. Math. 721 (2016) 109-147
\par
\medskip
\cite{Ch-N} S.G.Chiossi and P-A. Nagy {\it Complex homothetic
foliations on K\"ahler manifolds}  Bull. London Math. Soc. 44
(2012) 113-124.

\par
\medskip
\cite{D-M-1} A. Derdzi\'nski, G. Maschler {\it Special
K\"ahler-Ricci potentials on compact K\" ahler manifolds}, J.
reine angew. Math. 593 (2006), 73-116.
 \par
\medskip
\cite{D-M-2}  A. Derdzi\'nski, G.  Maschler {\it Local
classification of conformally-Einstein  K\"ahler metrics in higher
dimensions}, Proc. London Math. Soc. (3) 87 (2003), no. 3,
779-819.

\par
\medskip
\cite{J-1} W. Jelonek {\it K\"ahler manifolds with homothetic foliation by curves} Diff. Geom. and its Applications, no.46, (2016), 119-131.

\par
\medskip
\cite{J-2} W. Jelonek {\it K\"ahler manifolds with quasi-constant holomorphic curvature} Ann. Glob. Analysis and Geom., vol.36, (2009), 143-159.
par
\medskip
\cite{J-3} W. Jelonek {\it Killing tensors and warped products} Ann. Polon. Math., vol.LXXV.2, (2000), 111-124.

\cite{V} I. Vaisman {\it Conformal foliations} Kodai J. Math.
2(1979), 26-37.
\par
\medskip
Institute of Mathematics

Cracow University of Technology

Warszawska 24

31-155      Krak\'ow,  POLAND.

 E-mail address: wjelon\@pk.edu.pl
\bigskip

\enddocument